\documentclass[12pt,a4paper]{article}
\usepackage{amsthm}
\usepackage{amsfonts}
\usepackage{amssymb}
\usepackage{amsmath}
\usepackage{graphicx}
\usepackage{color}
\usepackage{colortbl}

\theoremstyle{definition}
\newtheorem{theorem}{Theorem}

\newtheorem{proposition}[theorem]{Proposition}
\newtheorem{lemma}[theorem]{Lemma}
\newtheorem{corollary}[theorem]{Corollary}

\newtheorem{conjecture}[theorem]{Conjecture}
\newtheorem{\proofname}{Proof}
\def\gcd{{\rm gcd}}

\title{\bf On the number of subsequences with a given sum in a finite
           abelian group}
\author{Gerard Jennhwa Chang,$^{123}$\thanks{E-mail: gjchang@math.ntu.edu.tw.
           Supported in part by the National Science Council under grant
           NSC98-2115-M-002-013-MY3.} ~
        Sheng-Hua Chen,$^{13}$\thanks{E-mail: b91201040@ntu.edu.tw.} \\
        Yongke Qu,$^4$\thanks{E-mail: quyongke@sohu.com.} ~
        Guoqing Wang,$^4$\thanks{E-mail: gqwang1979@yahoo.com.cn. Supported by NSFC (11001035).} ~
        and
        Haiyan Zhang$^5$\thanks{E-mail: yanhaizhang2222@sohu.com.}
        \\~ \\
{\scriptsize $^1$Department of Mathematics, National Taiwan University,
            Taipei 10617, Taiwan} \\[-0.15cm]
{\scriptsize $^2$Taida Institute for Mathematical Sciences, National
            Taiwan University, Taipei 10617, Taiwan} \\[-0.15cm]
{\scriptsize $^3$National Center for Theoretical Sciences, Taipei
            Office} \\[-0.15cm]
{\scriptsize $^4$Center for Combinatorics, LPMC-TJKLC, Nankai
University, Tianjin
            300071, P.R. China} \\[-0.15cm]
{\scriptsize $^5$Department of Mathematics, Harbin University of Science
            and Technology, Harbin 150080, P.R. China}}

\date{}

\begin{document}
\baselineskip 16pt
\parskip 5pt
\parindent 1truecm
\maketitle

\begin{abstract}
Suppose $G$ is a finite abelian group and $S$ is a sequence of
elements in $G$. For any element $g$ of $G$, let $N_g(S)$ denote the
number of subsequences of $S$ with sum $g$. The purpose of this
paper is to investigate the lower bound for $N_g(S)$.  In
particular, we prove that either $N_g(S)=0$ or $N_g(S) \ge
2^{|S|-D(G)+1}$, where $D(G)$ is the smallest positive integer
$\ell$ such that every sequence over $G$ of length at least $\ell$
has a nonempty zero-sum subsequence. We also characterize the
structures of the extremal sequences for which the equality holds
for some groups.
\end{abstract}

\newpage
\section{Introduction}              

Suppose $G$ is a finite abelian group and $S$ is a sequence over
$G$. The enumeration of subsequences with certain prescribed
properties is a classical topic in Combinatorial Number Theory going
back to Erd\H{o}s, Ginzburg and Ziv \cite{Erdos, Geroldinger1,
Geroldinger2} who proved that $2n-1$ is the smallest integer such
that every sequence $S$ over a cyclic group $C_n$ has a subsequence
of length $n$ with zero-sum. This raises the problem of determining
the smallest positive integer $\ell$ such that every sequence $S$ of
length at least $\ell$ has a nonempty zero-sum subsequence. Such an
integer $\ell$ is called the {\it Davenport constant}
\cite{Davenport} of $G$, denoted by $D(G)$, which is still unknown
in general.

For any $g$ of $G$, let $N_g(S)$ denote the number of subsequences
of $S$ with sum $g$. In 1969, J. E. Olson \cite{Olson} proved that
$N_0(S) \ge 2^{|S|-D(G)+1}$ for every sequence $S$ over $G$ of
length $|S| \ge D(G)$. Subsequently, several authors
\cite{Bialostocki, Balandraud, Cao, Furedi, Gao1, Gao2, Gao4, Gao5,
Grynkiewicz1, Grynkiewicz2, Guichard, Kisin} obtained a huge variety
of results on the number of subsequences with prescribed properties.
However, for any arbitrary $g$ of $G$, the lower bound of $N_g(S)$
remains undetermined.

In this paper, we determine the best possible lower bound of
$N_g(S)$ for an arbitrary $g$ of $G$. We also characterize the
structures of the extremal sequences which attain the lower bound
for some groups.

\section{Notation and lower bound}     

Our notation and terminology are consistent with \cite{Gao3}. We
briefly gather some notions and fix the notation concerning
sequences over abelian group. Let $\mathbb{N}$ and $\mathbb{N}_0$
denote the sets of positive integers and non-negative integers,
respectively. For integers $a,b \in \mathbb{N}_0$, we set $[a,b] =
\{x\in \mathbb{N}_0: a \le x \le b\}$. Throughout, all abelian
groups are written additively. For a positive integer $n$, let $C_n$
denote a cyclic group with $n$ elements.

Suppose $G$ is an additively written, finite abelian group and
$\mathcal{F}(G)$ is the multiplicatively written, free abelian
monoid with basis $G$. The elements of $\mathcal{F}(G)$ are called
{\it sequences} over $G$. We write sequences $S$ of $\mathcal{F}(G)$
in the form
$$
   S = \prod_{g\in G} g^{\upsilon_g(S)},
   \mbox{ with } \upsilon_g(S) \in \mathbb{N}_0 \mbox{ for all } g \in G.
$$
We call $\upsilon_g(S)$ the {\it multiplicity} of $g$ in $S$, and we
say that $S$ {\it contains} $g$ if $\upsilon_g(S) > 0$. A sequence
$T$ is called a {\it subsequence} of $S$ if $T|S$ in
$\mathcal{F}(G)$, or equivalently, $\upsilon_g(T) \le \upsilon_g(S)$
for all $g \in G$. The unit element $1 \in \mathcal{F}(G)$ is called
the {\it empty sequence}, and it is a zero-sum subsequence of any
sequence $S\in\mathcal{F}(G)$.  For sequences $S_1,S_2,  \ldots, S_n
\in \mathcal{F}(G)$, denote  by $\gcd(S_1,S_2,\ldots,S_n)$ the
longest subsequence dividing all $S_i$. If a sequence $S \in
\mathcal{F}(G)$ is written in the form $S = g_1g_2 \cdots  g_m$, we
tacitly assume that $m \in \mathbb{N}_0$ and $g_1,g_2,  \ldots, g_m
\in G$. For a sequence $S = g_1 g_2\cdots g_m = \prod_{g \in G}
g^{\upsilon_g(S)}$, we call
\begin{itemize} \itemsep 0pt
\item $|S| = m = \sum_{g \in G} \upsilon_g(S)$ the {\it length} of
$S$,

\item $\sigma(S) = \sum_{i=1}^{m} g_i$ the {\it sum} of $S$,

\item $\sum(S) = \{\sum_{i\in I} g_i: \phi \ne I \subseteq
[1,m]\}$ the {\it set of subsums} of $S$.
\end{itemize}
For convenience, we set $\sum^{\bullet}(S)=\sum(S)\cup \{0\}$ and
$-S = (-g_1)(-g_2)\cdots (-g_m)$. The sequence $S$ is called
\begin{itemize} \itemsep 0pt
\item {\it zero-sum free} if $0 \notin \sum(S)$,

\item a {\it zero-sum sequence} if $\sigma(S) = 0$,

\item a {\it minimal zero-sum sequence} if $S \ne 1$,
$\sigma(S)=0$, and every $T|S$ with $1 \le |T| < |S|$ is zero-sum
free.
\end{itemize}
For an element $g$ of $G$, let
$$
    N_g(S) = |\{I \subseteq [1,m]: \sum_{i\in I}g_i = g\}|
$$
denote the number of subsequences $T$ of $S$ with sum $\sigma(T)=g$,
counted the multiplicity of their appearance in $S$.  Notice that we
always have $N_0(S) \ge 1$.

\begin{theorem} \label{lower bound}
If $S$ is a sequence over a finite abelian group $G$ and $g \in
\sum^{\bullet}(S)$, then $N_g(S) \ge 2^{|S|-D(G)+1}$.
\end{theorem}

\noindent{\bf Proof.}  We shall prove the theorem by induction on
$m=|S|$. The case of $m \le D(G) -1$ is clear.  We now consider the
case of $m \ge D(G)$. Choose a subsequence $T|S$ of minimum length
with $\sigma(T)=g$, and a nonempty zero-sum subsequence
$W|T(-(ST^{-1}))$. By the minimality of $|T|$, $W$ is not a
subsequence of $T$, for otherwise $TW^{-1}$ is a shorter subsequence
of $S$ with $\sigma(TW^{-1})=g$. Choose a term $a|W$ but $a\nmid T$,
and let $X={\rm gcd}(W,T)$. Then, $-a|ST^{-1}$ such that $g =
\sigma(T) \in \sum^{\bullet}(S(-a)^{-1})$ and
$(g-\sigma(X))-(0-\sigma(X)-a)=g+a = \sigma(TX^{-1} (-(W(Xa)^{-1})))
\in \sum^{\bullet}(S(-a)^{-1})$. By the induction hypothesis,
$N_g(S) = N_g(S(-a)^{-1}) + N_{g+a}(S(-a)^{-1}) \ge
2^{m-D(G)}+2^{m-D(G)}=2^{m-D(G)+1}$. This completes the proof of the
theorem.    \qed

Notice that the result in \cite{Olson} that $N_0(S) \ge
2^{|S|-D(G)+1}$ for any sequence $S$ over $G$, together with the
following lemma, also gives Theorem \ref{lower bound}.
\begin{lemma} \label{transformation}
If $S$ is a sequence over a finite abelian group $G$, then for any
$T|S$ with $\sigma(T)=g \in \sum^{\bullet}(S)$,
  $$N_g(S)=N_0(T(-(ST^{-1}))).$$
\end{lemma}

\noindent{\bf Proof.} Let $\mathcal{A} = \{X|S: \sigma(X)=g\}$ and
$\mathcal{B} = \{Y|T(-(ST^{-1})): \sigma(Y)=0\}$. It is clear
that $|\mathcal{A}|=N_g(S)$ and
$|\mathcal{B}|=N_0(T(-(ST^{-1})))$. Define the map $\varphi:
\mathcal{A}\to \mathcal{B}$ by $\varphi(X)=TX_1^{-1}(-X_2)$ for
any $X \in \mathcal{A}$, where $X_1={\rm gcd}(X,T)$ and $X_2={\rm
gcd}(X,ST^{-1})$. It is straightforward to check that $\varphi$ is
a bijection, which implies $N_g(S)=N_0(T(-(ST^{-1})))$. \qed

We remark that the lower bound in Theorem \ref{lower bound} is best
possible. For any $g \in G$ and any $m \ge D(G)-1$, we construct the
extremal sequence $S\in \mathcal{F}(G)$ of length $m$ with
respect to $g$ as follows: Take a zero-sum free sequence $U \in
\mathcal{F}(G)$ with $|U|=D(G)-1$. Clearly, $U$ contains a
subsequence $T$ with $\sigma(T)=g$. For $S = T (-(UT^{-1}))
0^{m-D(G)+1}$, by Lemma \ref{transformation}, $N_g(S) = N_0(U
0^{m-D(G)+1}) = 2^{m-D(G)+1}$.

\begin{proposition} \label{one and all}
If $S$ is a sequence over a finite abelian group $G$ such that
$N_h(S) = 2^{|S|-D(G)+1}$ for some $h \in G$, then $N_g(S) \ge
2^{|S|-D(G)+1}$ for all $g\in G$.
\end{proposition}

\noindent{\bf Proof.} If there exists $g$ such that $N_g(S) <
2^{|S|-D(G)+1}$, then
  $$N_{h}(S(h-g))=N_{h}(S)+N_{g}(S)<2^{|S|+1-D(G)+1}$$
is a contradiction to Theorem \ref{lower bound} since $h \in
\sum^{\bullet}(S) \subseteq \sum^{\bullet}(S(h-g))$. \qed

\section{The structures of extremal sequences}    

In this section, we study sequence $S$ for which $N_g(G) =
2^{|S|-D(G)+1}$. By Lemma \ref{transformation},  we need only pay
attention to the case $g=0$.  Also, as $N_g(0S) = 2 N_g(S)$, it
suffices to consider the case $0 \nmid S$. For $|S|\geq D(G)-1$,
define $$E(S)=\{g\in G: N_g(S)=2^{|S|-D(G)+1}\}.$$

\begin{lemma} \label{ES structure}
Suppose $S$ is a sequence over a finite abelian group $G$ with
$0\nmid S$, $|S| \ge D(G)$ and $0 \in E(S)$. If $a$ is  a term of a
zero-sum subsequence $T$ of $S$, then
  $$E(S)+\{0, -a\}\subseteq E(Sa^{-1}).$$
\end{lemma}

\noindent{\bf Proof.} Since $0,-a \in \sum^{\bullet}(Sa^{-1})$,  by
Theorem \ref{lower bound}, $N_0(Sa^{-1}) \ge 2^{|S|-D(G)}$ and
$N_{-a}(Sa^{-1}) \ge 2^{|S|-D(G)}$. On the other hand,
$N_0(Sa^{-1})+N_{-a}(Sa^{-1})=N_0(S)=2^{|S|-D(G)+1}$ and so
$N_0(Sa^{-1})=N_{-a}(Sa^{-1})=2^{|S|-D(G)}$.  Hence, by Proposition
\ref{one and all}, $N_g(Sa^{-1}) \ge 2^{|S|-D(G)}$ for all $g\in G$.
Now, for every $h \in E(S)$, $N_h(Sa^{-1})+N_{h-a}(Sa^{-1}) = N_h(S)
=2^{|S|-D(G)+1}$ and so $N_h(Sa^{-1}) = N_{h-a}(Sa^{-1}) =
2^{|S|-D(G)}$, i.e., $ \{h, h-a\} \subseteq E(Sa^{-1})$. This proves
$E(S)+\{0, -a\} \subseteq E(Sa^{-1})$. \qed

\begin{lemma}\label{D(G)}(\cite{Geroldinger1}, Lemma 6.1.3, Lemma
6.1.4) Let $G \cong C_{n_1} \oplus C_{n_2} \oplus \cdots \oplus
C_{n_r}$ with $n_1|n_2|\cdots |n_r$, and $H$ be a subgroup of $G$,
then $D(G)\geq D(H)+D(G/H)-1$ and $D(G)\geq
\sum_{i=1}^{r}(n_i-1)+1$.
\end{lemma}

\begin{lemma} \label{nontrivial subgroup in Es}
If $S$ is a sequence over a finite abelian group $G$ such that
$E(S)$ contains a non-trivial subgroup $H$ of $G$, then $H \cong
\bigoplus_{i=1}^r C_2$ and $D(G)=D(G/H)+r.$
\end{lemma}

\noindent{\bf Proof.} Suppose $H \cong C_{n_1} \oplus C_{n_2} \oplus
\cdots \oplus C_{n_r}$, where $n_1 | n_2 | \cdots | n_r$, and assume
that $S=g_1g_2 \cdots g_m$. Consider the canonical map $\varphi: G
\to G/H$ and let $\varphi(S) = \varphi(g_1)\varphi(g_2) \cdots
\varphi(g_m)$ be a sequence over $G/H$. Then
$$
  |H| \cdot 2^{|S|-D(G)+1} = \sum_{h\in H} N_h(S) =
  N_0(\varphi(S)) \ge 2^{|\varphi(S)|-D(G/H)+1}.
$$
By  lemma \ref{D(G)}, we have $|H| \ge 2^{D(G)-D(G/H)} \ge
2^{D(H)-1}$, and so
$$
      \prod_{i=1}^r n_i \ge 2^{\sum_{i=1}^r (n_i-1)} = \prod_{i=1}^r 2^{n_i-1}.
$$
Hence, $n_i=2$ for all $i$, which gives $H\cong \bigoplus_{i=1}^r C_2$ and $D(G)=D(G/H)+r$.
\qed

Now, we consider the case $G=C_n$.  Notice that $D(C_n)=n$.

\begin{theorem} \label{C_n}
For $n \ge 3$, if $S\in \mathcal{F}(C_n)$ with $0\nmid S$ and
$N_0(S)=2^{|S|-n+1}$, then $n-1 \le |S|\leq n$ and  $S=a^{|S|}$, where
$a$ generates $C_n$.
\end{theorem}

\noindent{\bf Proof.} Suppose $S\in \mathcal{F}(C_n)$ with $0\nmid
S$ and $N_0(S)=2^{|S|-n+1}$. We first show  by induction that
\begin{equation}\label{supp(S)=a}
S=a^{|S|}
\end{equation} where $\langle a\rangle=C_n$. For $|S|=n-1$, we have
$N_0(S)=1$, i.e., $S$ is a zero-sum free sequence, and
\eqref{supp(S)=a} follows readily.

For $|S|\geq n$, since $N_0(S)=2^{|S|-n+1}\geq 2$, $S$ contains at
least one nonempty zero-sum subsequence $T$. Take an arbitrary
 term $c$ from $T$. By Lemma \ref{ES structure}, $0 \in E(Sc^{-1})$.
It follows from the induction hypothesis that $Sc^{-1}=a^{|S|-1}$
for some $a$ generating $C_n$. By the arbitrariness of $c$, we
conclude that \eqref{supp(S)=a} holds.

To prove $|S|\leq n$, we suppose to the contrary that $|S|\geq n+1$.
By $\eqref{supp(S)=a}$ and Lemma \ref{ES structure},
\begin{equation}\label{0 in E(a n+1)}
0\in E(a^{n+1}).
\end{equation} We see that $N_0(a^{n+1}) \geq
1+{n+1 \choose n}
> 4$, a contraction with \eqref{0 in E(a n+1)}.
\qed

Notice that Theorem \ref{C_n} is not true for $n=2$, since for any
$S\in \mathcal{F}(C_2)$ with $0\nmid S$, we always have $N_0(S)=2^{|S|-2+1}$.

While the structure of a sequence $S$ over a general finite abelian group $G$
with $0 \nmid S$ and $N_0(S)=2^{|S|-D(G)+1}$ is still not known, we have the
following result for the case when $|G|$ is odd.

\begin{theorem} \label{odd group}
If $S$ is a sequence over a finite abelian group $G$ of odd size
with $0\nmid S$ and $N_0(S)=2^{|S|-D(G)+1}$, then  $S$ contains
exactly $|S|-D(G)+1$ minimal zero-sum subsequences, all of which are
pairwise disjoint.
\end{theorem}

\noindent{\bf Proof.} We shall prove the theorem by induction on
$|S|$. If $|S|=D(G)$, then $N_0(S)=2$ and so $S$ contains exactly
one nonempty zero-sum subsequence. For $|S| \ge D(G)+1$, if all
minimal zero-sum subsequences of $S$ are pairwise disjoint, then the
conclusion follows readily. So we may assume that there exist two
distinct minimal zero-sum subsequences $T_1$ and $T_2$ with
$\gcd(T_1,T_2) \ne 1$. Take  a term $a|\gcd(T_1,T_2)$. By Lemma
\ref{ES structure}, $0 \in E(Sa^{-1})$ and so $Sa^{-1}$ contains
$r=|S|-D(G)\geq 1$ pairwise disjoint minimal zero-sum subsequences
$T_3, T_4, \ldots, T_{r+2}$ by the induction hypothesis.

We observe that no term $b$ is contained in exactly one $T_i$,
 where $i\in [1, r+2]$. Otherwise,   for
any $t\in [1, r+2]$, let $b|T_t$ and $b\nmid T_i$, where $i\in [1,
r+2]\setminus \{t\}$, by Lemma \ref{ES structure}, $0 \in
E(Sb^{-1})$ and then $Sb^{-1}$ contains  exactly $r$ minimal
zero-sum subsequences, which is a contradiction. Choose a term $c$
in $T_1$ but not in $T_2$.  By the above observation $c$ is in
another $T_i$, say $T_{r+2}$ and so not in any of $T_3, T_4, \ldots,
T_{r+1}$. Again $Sc^{-1}$ contains exactly $r$ disjoint minimal
zero-sum subsequences, which are just $T_2, T_3, \ldots, T_{r+1}$.
If $r \ge 2$, then  there exists a minimal zero-sum subsequence say
$T_{r+1}$, may only has common term with $T_1$.  By the observation,
we can get $T_{r+1}|T_1$, which is a contradiction to the
minimality. Thus $r=1$.

Now we have $N_0(S)$=4 and $T_1,T_2,T_3$ are the all minimal
zero-sum subsequences of $S$. If there is some
$d|\gcd(T_1,T_2,T_3)$, then $Sd^{-1}$ contains no minimal zero-sum
subsequence, which is impossible. Thus $\gcd(T_1,T_2,T_3) = 1$ and
$\sigma(\gcd(T_1,T_2)) = g \ne 0$, where $2g=0$, contradicting that
$|G|$ is odd. \qed

If we further assume that $E(S)=\{0\}$ in Theorem \ref{odd group},
the structure of $S$ can be  further restricted.

\begin{corollary}
If $S$ is a sequence over a finite abelian group $G$ of odd size
with $0 \nmid S$ and $E(S)=\{0\}$, then $S$ contains exactly
$r=|S|-D(G)+1$ minimal zero-sum subsequences, all of which are
pairwise disjoint, and nothing else, that is $S=T_1 T_2\cdots T_r$.
\end{corollary}

\noindent{\bf Proof.} By Theorem \ref{odd group},  $S$ contains
exactly $r=|S|-D(G)+1$ minimal zero-sum subsequences, all of which
are pairwise disjoint, that is $S=T_1T_2\cdots T_rW$. For any
subsequence $X$ of $S$ with $\sigma(X)=\sigma(W)$, if $W\nmid X$,
then $SX^{-1}$ is a zero-sum subsequence containing terms in $W$,
which is impossible.  So $W|X$, and then $\sigma(XW^{-1})=0$. This
gives $X=T_{i_1}T_{i_2}\cdots T_{i_s}W$ with $1\le i_1< i_2 < \cdots
<i_s \le r$. Hence, $N_{\sigma(W)}(S)=2^r$ and then $\sigma(W)\in
E(S) =\{0\}$ implying $W=1$. \qed

The property that $S$ contains exactly $|S|-D(G)+1$ minimal zero-sum
subsequences, all of which are pairwise disjoint, implies that $|S|$
is bounded as in the case of Theorem \ref{C_n} for cyclic groups.

\begin{theorem} \label{infinite sequence}
For any finite abelian group $G \cong C_{n_1} \oplus C_{n_2} \oplus
\cdots \oplus C_{n_r}$ with $n_1|n_2|\cdots |n_r$, (i) implies the
three equivalent statements (ii), (iii) and (iv).
\begin{enumerate} \itemsep 0pt
\item[(i)]   Any sequence $S$ over $G$ with $0\nmid S$ and $N_0(S)=2^{|S|-D(G)+1}$,
              contains exactly $|S|-D(G)+1$ minimal zero-sum subsequences, all of which are pairwise disjoint.

\item[(ii)]  There is a natural number $t=t(G)$ such that $|S|\le t$ for every
             $S\in \mathcal{F}(G)$ with $0\nmid S$ and $N_0(S)=2^{|S|-D(G)+1}$.

\item[(iii)] For any subgroup $H$ of $G$ isomorphic to $C_2$,
             $D(G)\geq D(G/H)+2$.
\item[(iv)]  For any $S\in \mathcal{F}(G)$,
             $E(S)$ contains no non-trivial subgroup of $G$.

\end{enumerate}
\end{theorem}

\noindent{\bf Proof.}
(i) $\Rightarrow$ (ii). Since $S$ contains exactly $|S|-D(G)+1$ minimal
zero-sum subsequences, all of which are pairwise disjoint, $|S| \ge 2(|S|-D(G)+1)$
which gives that $|S| \le 2D(G)-2$.

(ii) $\Rightarrow$ (iii).  Assume to the contrary that
$D(G)=D(G/H)+1$ for some subgroup $H=\{0,h\}$ of $G$. Let
 $\varphi: G
\to G/H$ be the canonical map, and $m=D(G/H)$, we choose a sequence
$S=g_1g_2\cdots g_m$ over $G$ such that $\varphi(S) =
\varphi(g_1)\varphi(g_2) \cdots \varphi(g_m)$ is a minimal zero-sum
sequence over $G/H$, and $\sigma(S)=h$ in $G$. Since
$$
  N_0(S)+N_h(S)=N_0(\varphi(S))=2=2 \cdot 2^{|S|-D(G)+1}
$$
and $N_0(S)$ and $N_h(S)$ are not zero,  by theorem \ref{lower
bound},  $N_0(S)=N_h(S)=2^{|S|-D(G)+1}$. Since $N_0(Sh^k) =
N_0(Sh^{k-1}) + N_h(Sh^{k-1}) = N_h(Sh^k)$, by induction we have
$N_0(Sh^k) = N_h(Sh^k) = 2^{|Sh^{k}|-D(G)+1}$ for all $k$, a
contradiction to the assumption in (ii).

(iii) $\Rightarrow$ (iv). Suppose to the contrary that there exists
a sequence $S\in \mathcal{F}(G)$ such that $E(S)$ contains a
non-trivial subgroup $H$ of $G$. By Lemma \ref{nontrivial subgroup
in Es}, $H \cong \bigoplus_{i=1}^s C_2$ and $D(G)=D(G/H)+s$. Hence,
$E(S)$ contains a subgroup $H' \cong C_2$.  If $D(G)\geq D(G/H')+2$,
then by lemma \ref{D(G)}, $D(G)\geq D(G/H')+2\geq
D(H/H')+D((G/H')/(H/H'))+1=s+1+D(G/H)>D(G)$, a contradiction.

(iv) $\Rightarrow$ (ii). For $|S|\geq D(G)$, that is,
$N_0(S)=2^{|S|-D(G)+1}>1$, there exists a nonempty zero-sum
subsequence $T_1$ of $S$ and a term $a_1|T_1$. By Lemma \ref{ES
structure}, $0\in E(S) \subseteq E(Sa_1^{-1})$. By (iv), $\langle -a_1
\rangle \not\subseteq E(Sa_1^{-1})$. Let $k$ be the minimum index
such that $k(-a_1) \notin E(Sa_1^{-1})$, that is,
$\{0,-a_1,\ldots,(k-1)(-a_1)\}\subseteq E(Sa_1^{-1})$ but
$k(-a_1)\notin E(Sa_1^{-1})$. Then,
$N_{(k-1)(-a_1)}(Sa_1^{-1})=2^{|Sa_1^{-1}|-D(G)+1}$ but
$N_{k(-a_1)}(Sa_1^{-1}) \ne 2^{|Sa_1^{-1}|-D(G)+1}$. Thus,
$$N_{(k-1)(-a_1)}(S)=N_{(k-1)(-a_1)}(Sa_1^{-1})+N_{k(-a_1)}(Sa_1^{-1})
\ne 2^{|S|-D(G)+1}$$ and so $(k-1)(-a_1)\notin E(S)$. This means
$$E(S)\subsetneq E(Sa_1^{-1}).$$ If $|Sa_1^{-1}| \geq D(G)$, a similar
argument shows that there exists a nonempty zero-sum
subsequence $T_2$ of $Sa_1^{-1}$ and a term $a_2|T_2$, thus,
$E(Sa_1^{-1})\subsetneq E(Sa_1^{-1}a_2^{-1})$. We continue this
process to get $a_1,a_2,\ldots,a_{|S|-D(G)+1}$ of $S$ such that
$$E(S)\subsetneq E(Sa_1^{-1})\subsetneq\cdots \subsetneq
E(Sa_1^{-1}a_2^{-1}\cdots a_{|S|-D(G)+1}^{-1}).$$ Since
$|E(Sa_1^{-1}a_2^{-1}\cdots a_{|S|-D(G)+1}^{-1})|\leq |G|$, we
conclude $|S|\leq D(G)+|G|-1:=t$. \qed

\section{Concluding remarks}        

We are interested in the structure of a sequence $S$ over a finite
abelian group $G$ such that $N_0(S) = 2^{|S|-D(G)+1}$. Based on the
experiences in Section 3, we have the following two conjectures.

\begin{conjecture}
Suppose $G$ is a finite abelian group in which $D(G) \geq D(G/H)+2$
for every subgroup $H$ of $G$ isomorphic to $C_2$.  If $S \in
\mathcal{F}(G)$  is a sequence with $0\nmid S$ and
$N_0(S)=2^{|S|-D(G)+1}$,
 then $S$ contains exactly $|S|-D(G)+1$ minimal zero-sum subsequences, all of
 which are pairwise disjoint.
\end{conjecture}

Notice that this conjecture holds when $G$ is cyclic or $|G|$ is odd.
The second conjecture concerns the length of $S$.

\begin{conjecture} \label{length-conj}
Suppose $G\cong C_{n_1}\oplus C_{n_2}\oplus \cdots \oplus C_{n_r}$
where $1<n_1|n_2|\cdots |n_r$ and $D(G)=d^*(G)+1=\sum_{i=1}^r
(n_i-1)+1$. If $S\in \mathcal{F}(G)$ with $0\nmid S$ and $E(S)\neq
\emptyset$ contains no non-trivial subgroup of $G$, then $|S|\leq
d^*(G)+r$.
\end{conjecture}

The following example shows that if Conjecture \ref{length-conj}
holds, then the upper bound $d^*(G)+r=\sum_{i=1}^r n_i$ is best
possible. Let $G\cong C_{n_1}\oplus C_{n_2}\oplus \cdots \oplus
C_{n_r}=\langle e_1\rangle \oplus \langle e_2\rangle \oplus \cdots
\oplus \langle e_r\rangle$ with $1<n_1|n_2|\cdots |n_r$. Clearly,
$S=\prod_{i=1}^r e_i^{n_i}$ is an extremal sequence with respect to
$0$ and of length $d^*(G)+r$.


\end{document}